\documentclass[10pt]{article}
\usepackage{amsmath, amssymb}
\usepackage{latexsym}
\usepackage{mathtools}
\usepackage{ulem}
\usepackage{graphicx}
\usepackage{color}
\usepackage{url}
\numberwithin{equation}{section}
\DeclareMathAlphabet{\itbf}{OML}{cmm}{b}{it}

\newcommand{\RR}{\mathbb{R}}

\newcommand{\ds}{\displaystyle}
\newcommand{\no}{\nonumber}
\newcommand{\ri}{\rightarrow}

\def\finproof{\hfill {\small $\Box$} \\}

\newcommand{\q}{\quad}

\newcommand{\bt}{\mbox{\boldmath$\tau$}}

\newcommand{\bp}{{\itbf p}}
\newcommand{\bq}{{\itbf q}}

\newcommand{\bx}{{\itbf x}}
\newcommand{\bv}{{\itbf v}}
\newcommand{\bg}{{\itbf g}}
\newcommand{\bh}{{\itbf h}}
\newcommand{\bev}{{\itbf e}}
\newcommand{\bw}{{\itbf w}}
\newcommand{\bu}{{\itbf u}}
\newcommand{\bn}{{\itbf n}}
\newcommand{\by}{{\itbf y}}

\newcommand{\bi}{\begin{itemize}}
\newcommand{\ei}{\end{itemize}}

\newcommand{\bH}{{\itbf H}}

\newcommand{\bG}{\itbf G}
\newcommand{\bK}{\itbf K}

\newcommand{\be}{\begin{eqnarray}}
\newcommand{\ee}{\end{eqnarray}}

\newcommand{\ben}{\begin{eqnarray*}}
\newcommand{\een}{\end{eqnarray*}}
\def\ds{\displaystyle}

\newcommand\ov{\overline}

\newtheorem{lem}{Lemma}[section]

\newtheorem{thm}{Theorem}[section]

\newcommand{\bea}{\begin{eqnarray*}}
\newcommand{\eea}{\end{eqnarray*}}
\newcommand{\bean}{\begin{eqnarray}}
\newcommand{\eean}{\end{eqnarray}}
\newcommand{\p}{\partial}
\newcommand{\f}{\frac}
\newcommand{\s}{\sqrt}
\newcommand{\di}{\mbox{div }}
\newcommand{\aaa}{\mbox{$[$}}
\newcommand{\bbb}{\mbox{$]$}}

\begin{document}

\title{Stability estimates for the fault inverse problem}

\author{ Faouzi Triki,\thanks{Laboratoire Jean Kuntzmann,  UMR CNRS 5224, 
Universit\'e  Grenoble-Alpes, 700 Avenue Centrale,
38401 Saint-Martin-d'H\`eres, France. Email: faouzi.triki@univ-grenoble-alpes.fr.} 
 \and Darko
Volkov \thanks{Department of Mathematical Sciences,
Worcester Polytechnic Institute, Worcester, MA 01609. Email: darko@wpi.edu.
}  }

\maketitle
\tableofcontents 

\begin{abstract}
We study in this paper stability estimates for the fault inverse problem. 
In this problem, faults are assumed to be  planar open surfaces in a half space elastic medium
with known Lam\'e coefficients. A traction free condition is imposed 
on the boundary of the half space. 
Displacement fields present jumps across  faults, 
called  slips,
while  traction derivatives are continuous.
It was proved in \cite{volkov2017reconstruction} that 
if the displacement field is known on an open set on the boundary of the half space, then the fault 
and the slip are uniquely determined. 
In this present paper, we study the stability of this uniqueness result with regard 
to the coefficients of the equation of the plane containing the fault. 
If the slip field is known we  state and prove a Lipschitz stability result.
In the more interesting  case where the slip field is unknown,
 we  state and prove another Lipschitz stability result 
under the additional assumption, which is still physically relevant,
that the slip field is one directional.
\end{abstract}

\bigskip



\section{Introduction}
Understanding and mapping the structure of Earth's crust 
in ever finer details has always captured the interest of geophysicists.
Seismic and displacements data  are collected by sensors and then processed 
using  Partial Differential Equations (PDE)  models and inverse problem formulations. Typical 
 models for the Earth's crust 
involve linear elasticity equations: this is because displacements and deformations
are very small compared to the thickness of the crust. Moreover, if local phenomena such as
earthquakes or active subduction zones are studied,  a half space formulation is adequate
\cite{arnadottir1994, fukahata2008, jonsson2002fault}.
With the advent of ultra accurate 
satellite based measurements of surface displacements (2 to 5 millimeter resolution)
the study of so called "slow earthquakes"
\cite{Cavalie2013, Dragert2001,  GuilhemNadeau2012, fu2013repeated,  
Radiguet, Radiguet2011} has recently attracted a lot of attention.
Most   authors  first set   a profile for the interface between tectonic plates (also called faults)
derived from seismicity or gravimetry as in \cite{Kim} and  then use a linear inverse algorithm
for determining   slip fields on faults. 
A popular algorithm is the one explained in Tarantola's textbook \cite{tarantola2005inverse}.
In addition to recovering these slip fields from surface displacement measurements, some authors 
have sought to simultaneously recover some geometric features of the fault, such as the dip angle
\cite{arnadottir1994, fukahata2008}.
However, until recently, there was no formal mathematical proof that the 
simultaneous recovery of the (piecewise linear) geometry of the 
fault and the slip was at all possible.
This was achieved in \cite{volkov2017reconstruction}. 
From there the second author and Sandmunienge have derived a deterministic and a stochastic
fault reconstruction algorithm
\cite{volkov2017stochastic} and estimated 
convergence rates to the  solution of the inverse problem.
In \cite{volkov2017stochastic}, these
convergence rates
 still depend on  the intrinsic  stability of the underlying inverse problem.
Although numerical investigations hinted at a possible Lipschitz type stability, these stability estimates
were still unknown at the time of writing of \cite{volkov2017stochastic}, and 
they are the subject of this present
study. In general,
a uniqueness statement for solving an inverse problem is not 
of great practical use without a stability result. 
From a pragmatic and computational point of view, mathematical objects can only 
be computed approximately and real life field data is always tinted by measurement errors, so one would not
want these errors to grow exponentially in inversion algorithms.\\
A literature review of the field of 
  inverse problems 
will show that stability results 
 are  notoriously difficult to derive and prove.
The major difficulty in proving such results is that   solutions to  inverse problems are not explicitly formulated.
There are a few papers on stability estimates for the recovery of cracks in materials,
which is the analog of faults in Earth's crust.
In an earlier paper,  Friedman and Vogelius  \cite{friedman1989determining}
showed a stability result for the recovery of liner cracks. 
In that paper the governing equation was the two dimensional conductivity and
outer boundary conditions were prescribed to adequate values.
In \cite{alessandrini2000classe}, Alessandrini et al. proved a general $\log \log$ stability 
estimate for the Hausdorff distance between two $C^{1,1}$ domains where
in each domain there is a solution to the same conductivity equation with same Neumann condition
and the stability estimate is in the $L^2$ distance between the corresponding Dirichlet outputs
on one part of the boundary. In \cite{ammari2017identification}, Ammari and the first author 
were able to  improve this $\log \log$ estimate based on the assumptions
that Dirichet data is available for a whole
range of frequencies and boundaries are a priori known to be open real analytic curves.
In \cite{beretta2008determination}, Beretta et al. were also able to derive and prove an
interesting Lipschitz stability result.  Their result pertains to two dimensional linear elasticity
in bounded domains with cracks. 
In the case of linear cracks they were able to derive Lipschitz continuity of the 
Hausdorff distance between cracks in terms of overdetermined boundary data. Here we need
to explain that the case of faults which pertains to our research project is drastically
different  since we can not impose boundary conditions, we are simply
passively measuring displacements on one part of the boundary of an infinite domain
while an unknown slip field on an unknown fault is the forcing term of our governing
equations.\\\\
This paper is organized as follows. In section \ref{Mathematical} we introduce the PDE for the forward
fault problem and we recall the uniqueness statement for the inverse fault problem proved in
\cite{volkov2017reconstruction}. Section 
\ref{first stability} contains  our first stability result and its proof. 
It relies on the implicit function theorem.
Indeed, using the Green's function for the forward PDE, solutions can represented
by convolution with the slip on the fault. We thus define a fault to surface operator. 
At fixed slip, in effect, this introduces
 a $C^1$ function $\phi$ from the set of geometry parameters $m$
in $\RR^3$ defining the plane containing the fault to the space of surface measurements.
We know from   \cite{volkov2017reconstruction} that $\phi$ is injective. 
Thus, by the inverse function Theorem, if $\nabla \phi (m)$ has   full rank, 
$\phi^{-1} $ is $C^1$ in a neighborhood of $\phi(m)$ and is therefore Lipschitz continuous.
The crux of the proof is in proving that $\nabla \phi (m)$ has full rank. 
This is established by an argument by contradiction. If $\nabla \phi (m)$ does not have full rank, then 
 using relations on jumps for the Green's function of the forward problem (and of its derivatives),
we can derive a PDE for the slip field $\bh$ on the fault. Finally, we prove
that this PDE can only have the trivial solution, completing the proof. 
It turns out that although the PDE on $\bh$ is a relatively simple transport equation in 
most cases,
in the particular case of horizontal faults a much more complicated system of PDE for $\bh$ must be solved.
In section \ref{Second stability theorem}, we assume that a 
fixed but 
\textsl{unknown}  slip $\bh_0$  is occurring on a plane with geometry parameter $m_0$.
Thus in this case it is not possible to evaluate the difference $ \|\phi(m) - \phi(m_0)\|$.
Instead, we use a linear operator $A_m$ mapping any slip $\phi$ to surface measurements 
$A_m \phi$ and we may minimize
$\inf \| A_m \phi - A_{m_0} \bh_0 \|$ where the $\inf$ is taken over all possible slips.
This quantity is proven to be bounded below by a constant times $|m -m_0|$ under the additional assumption
 that $\bh_0 $ is one directional (that direction is not known for the inverse problem), or 
that $\bh_0 $ is a gradient. Finally, the rather technical formulas for the jumps
of integrals containing a convolution of  the elasticity Green's tensor with a vector field density
are shown in Appendix. Although a related formula is  standard in  solid mechanics,
we have not found formulas for the jumps of the first and second derivatives in the literature. 
This is probably due to the fact that they are not directly related to a physical problem 
and that they may be too intricate to prove without the use of symbolic computation software.

\section{Mathematical model and uniqueness result} \label{Mathematical}

\subsection{Forward problem}
Using the standard rectangular coordinates $\bx= (x_1, x_2, x_3)$ of $\RR^{3}$,
we define $\RR^{3-}$ to be the open half space $x_3<0$.
Let $\p_i$ denote the derivative in the $i$-th coordinate.
In this paper 
we  consider the case of linear, homogeneous, isotropic elasticity; 
the two Lam\'e constants $\lambda$ and $\mu$ will be two positive constants.
For a vector field $\bu = (u_1, u_2, u_3)$
the stress  and strain tensors will be denoted as follows,
\bea
\sigma_{ij}(\bu) = \lambda \, \di \bu \, \delta_{ij} + \mu \, (\p_i u_j + \p_j u_i ), \\
\epsilon_{ij}(\bu) = \f12 (\p_i u_j + \p_j u_i ),
\eea
and the  stress vector in the 
direction $\bev \in \mathbb R^3 $ will be denoted by
\bea
T_\bev \bu = \sigma (\bu) \bev.
\eea

Let $\Gamma$ be a Lipschitz open surface which is strictly included in $\RR^{3-}$, with a normal
vector $\bn$.
We define the jump  $\aaa \bv \bbb$ of a vector field $\bv$ across $\Gamma$ 
to be
$$
\aaa \bv \bbb (\bx) = \lim_{h \ri 0^+} \bv  (\bx + h \bn) - \bv  (\bx - h \bn),
$$
for $\bx$ in $\Gamma$, if this limit exists.
Let $\bu$ be the displacement field solving
\bean
\mu \Delta \bu+ (\lambda+\mu) \nabla \di \bu= 0  \mbox{ in } \RR^{3-} 
\setminus \Gamma \label{uj1}, \\
 T_{\bev_3} \bu =0 \mbox{ on the surface } x_3=0 \label{uj2}, \\
 T_{\bn} \bu  \mbox{ is continuous across } \Gamma \label{uj3}, \\
 \aaa  \bu\bbb =\bg \mbox{ is a given jump across } \Gamma , \label{uj3andhalf}\\
\bu (\bx) = O(\f{1}{|\bx|^2}),  \nabla \bu (\bx) = O(\f{1}{|\bx|^3}), \mbox{ uniformly as } 
|\bx| \rightarrow \infty,
\label{uj4}
\eean
where 
$\bev_3$ is the vector $(0,0,1)$.
	For vector fields $\bv,\bw$ in $\RR^{3 -} \setminus \ov{\Gamma}$  
whose gradient is square integrable we introduce the bilinear product
\bea
B(\bv,\bw)=\int_{\RR^{3-}\setminus \ov{\Gamma}} \lambda\,  \mbox{tr}(\nabla \bv)  \mbox{tr}(\nabla \bw)  
 + 2 \mu \, \mbox{tr} ( \epsilon(\bv) \epsilon(\bw) ),
\eea
where $\mbox{tr} $ is the trace. 
In \cite{volkov2017reconstruction}, we defined the functional space 
$\cal{V}$ of vector fields $\bv $ defined in $\RR^{3-}\setminus{\ov{\Gamma}}$
such that $\nabla \bv $ and $\ds \f{\bv}{(1+ r^2)^{\f12}}$ are in $L^2(\RR^{3-}\setminus{\ov{\Gamma}})$, and we proved that 
the following four norms are equivalent on $\cal{V}$:
\bea
\| \bv \|_1 = (\int_{\RR^{3-}\setminus{\ov{\Gamma}}} |\nabla \bv|^2)^{\f12} +
 (\int_{\RR^{3-}\setminus{\ov{\Gamma}}} \f{| \bv|^2}{1 + r^2})^{\f12}, \\
\| \bv \|_2 = (\int_{\RR^{3-}\setminus{\ov{\Gamma}}} |\nabla \bv|^2)^{\f12}, \\
\| \bv \|_3 = (\int_{\RR^{3-}\setminus{\ov{\Gamma}}} |\epsilon ( \bv)|^2)^{\f12}, \\
\| \bv \|_4 = B(\bv, \bv)^{1/2}. \\
\eea
 Let  $D$ be a bounded domain with a  Lipschitz boundary $\p D$ containing $\Gamma$. We define the Sobolev space 
  $\widetilde{H}^{\f12}(\Gamma)^2$ 
 to be the set of restrictions to $\Gamma$ of tangential fields  in
   $H^{\f12}(\p D)^2$ supported in $\Gamma$.
We  proved in  \cite{volkov2017reconstruction} the following theorem.
\begin{thm} Let $\bg$ be in $\widetilde{H}^{\f12}(\Gamma)^2$.
The problem
(\ref{uj1}-\ref{uj3andhalf})
 has a unique solution in ${\cal V}$.
In addition the solution $\bu$ satisfies the decay conditions
(\ref{uj4}).
\end{thm}

In this paper we will only consider forcing terms $\bg$ which are tangential to 
$\Gamma$. Physically, this suggests that the fault $\Gamma$ is not opening or starting to self intersect: only slip is allowed.
 We recall that if $\bg$ is continuous, the support of $\bg$, 
$\mbox{supp } \bg$, is equal to the closure
of the set of points in $\Gamma$ where $\bg$ is non zero; in general $\mbox{supp }  \bg$
is defined in the sense of distributions. 

\subsection{Fault inverse problem} \label{Fault inverse problem}
Can we determine both $\bg$ and $\Gamma$ from the data
 $\bu$ given only on the plane $x_3=0$?
The following Theorem in \cite{volkov2017reconstruction}
asserts that this is possible if the data is known on a relatively open
set  of the plane $x_3=0$.
\begin{thm}{\label{uniq1}}
Let $\Gamma_1$ and $\Gamma_2$ be two bounded open 
surfaces, with smooth boundary, such that each of them  is included in a
rectangle strictly contained in $\RR^{3-}$.
For $i$ in $\{ 1,  2\}$, assume that $\bu^i$ solves  (\ref{uj1}-\ref{uj4}) for 
$\Gamma_i$ in place of $\Gamma$ and $\bg^i$, a tangential field in 
$H^{1}_0(\Gamma_i)^2$, in place of $\bg$.
Assume that $\bg^i$ has full support in $\Gamma_i$, that is, 
$\mbox{supp } \bg_i = \ov{\Gamma_i}$.
Let $V$ be a non empty open subset in $\{x_3 =0\}$.
If $\bu^1$ and $\bu^2$ are equal in $V$, then
$\Gamma_1 =  \Gamma_2$
and $\bg^1=\bg^2$.
\end{thm}
The  solution $\bu$ to problem (\ref{uj1}-\ref{uj3andhalf})  can also be written out 
 as the convolution on $\Gamma$
\bean
  \int_\Gamma \bH(\bx, \by, \bn) \bg(\by) \, d \sigma (\by) \label{int formula},
\eean
where $\bH$ is  the  Green's tensor  associated to the system (\ref{uj1}-\ref{uj4}), and
  $\bn$ is the normal to $\Gamma$. 
The practical determination of this adequate half space Green's tensor $\bH$ was first studied
in  \cite{Okada} and later, more rigorously, in \cite{DV}.
In particular, $\bH$ satisfies the decay conditions
\bea
 \bH(\bx, \by, \bn) = O(|\bx|^{-2}), \, \nabla_x  \bH(\bx, \by, \bn) = O(|\bx|^{-3}), \q |\bx| \ri \infty,
\eea
uniformly in $\by$ and in $\bn$, as long as $\by$ remains in a bounded subset of
$\RR^{3-}$.
Due to formula (\ref{int formula}) 
we can define a continuous mapping ${\cal M}$ from tangential fields
$\bg $ in $H^{1}_0(\Gamma)^2$ to surface displacement fields $\bu(x_1, x_2, 0)$
in $L^2(V) $ where $\bu$ and $\bg$ are related
by (\ref{uj1}-\ref{uj4}). Theorem \ref{uniq1} asserts that 
this mapping is injective, so an inverse operator can be defined. 
It is well known, however, that such an operator ${\cal M}$ is compact, therefore
its inverse is unbounded. It is thus clear that any stable numerical method for reconstructing
$\bg$ from $\bu(x_1, x_2, 0)$ will have to use some regularization process.
Our goal in this paper is to analyze the stability properties of the fault inverse problem
with regard to the plane containing $\Gamma,$ first in section 
\ref{first stability} as the slip on the fault is fixed, and then in 
section \ref{Second stability theorem} in the case of unknown slips.
 
\section{Lipschitz stability of the fault geometry for a fixed slip}\label{first stability}

\subsection{Preliminary results}
The formula for the Green tensor $\bH(\bx, \by, \bn)$ (\ref{int formula}) is given in 
\cite{Okada, DV}. 
Here we only give an explicit formula for its free space analog $\bG (\bx, \by, \bn)$
and we will use the fact that the difference 
$\bH(\bx, \by, \bn) -\bG (\bx, \by, \bn)$ is a smooth function for $(\bx, \by)$
in $\RR^{3-} \times \RR^{3-} $.
Recall the well known formula for Kelvin's Green's tensor
\bean
\bK_{ij}(\bx,\by)=\f{1}{8\pi \mu (\lambda+2 \mu)}
( (\lambda + \mu ) \p_{x_i }\, r \p_{x_j } r + (\lambda + 3 \mu ) \delta_{ij} )
\f{1}{r}, \label{kelvin}
\eean
where $r=|\bx -\by|$. Now if $\bv$ is any fixed vector, define
\bean \label{Gdef}
\bG(\bx,\by,\bv)=(T_{\bv(y)} \bK(\bx,\by))^T .
\eean 
We will need the following formulas for the jumps across $\Gamma$ of vector fields defined
by surface convolution of densities against $\bG$.
\begin{lem} \label{all jumps}
Let $\Gamma$ be an open surface in $\RR^3$ 
included in the plane $x_3=0$.
Let $\bg$ be  a three dimensional 
vector field on $\Gamma$ with regularity  $C^\infty_c(\Gamma)$. 
Then the following jump formulas across 
$\Gamma$
hold
\bean
 \aaa  \int_\Gamma \bG(\bx, \by, \bev_3) \bg(\by) \, d y_1 dy_2 \bbb = \bg(\bx), 
\label{jump1}\\
\aaa  \int_\Gamma (\p_{y_1}\bG)(\bx, \by, \bev_3) \bg(\by) \, d y_1 dy_2 \bbb = -  
\p_{x_1}\bg(\bx), 
\label{jump2}\\
\aaa  \int_\Gamma \bG(\bx, \by, \bev_1) \bg(\by) \, d y_1 dy_2 \bbb = 
\f{\lambda}{\lambda  + 2 \mu} g_1 (\bx)  \bev_3
 + g_3 (\bx) \bev_1, \label{jump3}\\
\aaa  \int_\Gamma (\p_{y_3}\bG)(\bx, \by, \bev_3) \bg(\by) \, d y_1 dy_2\bbb = 
\f{\lambda}{\lambda  + 2 \mu} (\di_{\Gamma}\bg_\Gamma ) (\bx ) \bev_3 + 
\nabla_\Gamma g_3 (\bx),\label{jump4} \\
\aaa  \p_{x_3} \int_\Gamma  \bG(\bx, \by, \bev_3) \bg(\by) \, d y_1 dy_2\bbb = 
-\f{\lambda}{\lambda  + 2 \mu} (\di_{\Gamma}\bg_\Gamma ) (\bx ) \bev_3 - 
\nabla_\Gamma g_3 (\bx),\label{jump4+1/2} 
\eean
where $\bg = (g_1, g_2,g_3)$, $\bg_\Gamma = (g_1, g_2, 0)$, 
$\di_{\Gamma}\bg_\Gamma = \p_1 g_1 + \p_2 g_2$, and $\nabla_\Gamma ( \bg \cdot \bev_3)
= (\p_1 g_3, \p_2 g_3, 0 )$.
For the  normal derivative of (\ref{jump4}) we have the jump formula
\bean
\aaa \p_{x_3} \int_\Gamma (\p_{y_3} \bG)(\bx, \by, \bev_3) \bg(\by) \, d y_1 dy_2 \bbb = 
(\f{3\lambda + 4 \mu}{\lambda + 2 \mu} \p_{1}^2 g_1 + \p_{2}^2 g_1 + 2 
\f{\lambda + \mu}{\lambda + 2 \mu} \p_{1}\p_{2} g_2) \bev_1  \no \\
+ (\f{3\lambda + 4 \mu}{\lambda + 2 \mu} \p_{2}^2 g_2 + \p_{1}^2 g_2 + 2 
\f{\lambda + \mu}{\lambda + 2 \mu} \p_{1}\p_{2} g_1) \bev_2 \no \\
-\f{\lambda }{\lambda + 2 \mu} (\p_1^2 + \p_2^2) g_3 \bev_3.
\label{jump5}
\eean
Finally, we give a jump formula for the normal derivative of
(\ref{jump3}) 
\bean
\aaa \p_{x_3} \int_\Gamma \bG(\bx, \by, \bev_1) \bg(\by) \, d y_1 dy_2 \bbb &= &
\label{jump6} (\f{3 \lambda + 4 \mu}{\lambda + 2 \mu} \p_1 g_1 + \p_2 g_2 ) \bev_1
+(\f{\lambda}{\lambda + 2 \mu } \p_2 g_1 + \p_1 g_2) \bev_2 \no \\
&&- \f{\lambda}{\lambda + 2 \mu }  \p_1 g_3 \bev_3.
\eean
\end{lem}
\textbf{Proof:}
Formula (\ref{jump1}) is well known, however, formulas (\ref{jump2}-\ref{jump6}) are not
readily found in the literature. It is therefore worth providing a proof, which can be 
found in the Appendix.

\begin{lem} \label{jump lem2}
Let $\bg$ be be a tangential vector field on $\Gamma$ with regularity  $H^1_0(\Gamma)$. 
Then the jump formulas (\ref{jump1}) and (\ref{jump3}) still hold in the $H^1_0(\Gamma)$ norm,
while the jump formulas  (\ref{jump2}, \ref{jump4}, \ref{jump4+1/2}, \ref{jump6}) hold 
as  continuous linear operations
 from $H^1_0(\Gamma)$
to $L^2(\Gamma)$. The jump formula (\ref{jump5}) holds as a continuous linear
operator from $H^1_0(\Gamma)$ to $H^{-1}(\Gamma)$.
\end{lem}
\textbf{Proof:} This is clear since $C^\infty_c(\Gamma)$ is dense in $H^1_0(\Gamma)$.

\begin{lem} \label{deg 1 eqn lem}
Let $\Omega$ be an open bounded subset in $\RR^2$. Let $f$ be in $C^\infty(\RR^2)$
such that $f^{-1}(\{ 0 \}) \cap \Omega$ has zero measure. Let $\bt$ be a non-zero vector 
in $\RR^2$. Assume that $u$ is in $H^1_0(\Omega)$ and satisfies in $\Omega$ the partial
differential equation
\bean \label{deg 1 eqn}
 \p_{\bt} ( f u  )
+ \alpha u = 0,
\eean
where $\alpha$ is a constant in $\RR$.
Then $u$ is zero.
\end{lem}
\textbf{Proof:}
We first assume that $\alpha \neq 0$. 
We  note that for any function $g$ in $C^\infty_c(\Omega)$ 
the divergence theorem implies that
\bean \label{div th}
0 = \int_\Omega \di ( g \bt) = \int_\Omega \nabla g \cdot \bt,
\eean
which can be extended by density to any $g$ in $H^1_0(\Omega)$.
Let $f_n$ be a sequence in $C^\infty(\Omega)$ which converges to $f^+$ in the
$H^1$ norm. 
By formula (\ref{div th}),
\bean \label{div th2}
  0=  \int_{\Omega} \nabla (f_n f u^2) \cdot \bt =
	\int_{\Omega}  f_n u \nabla ( f u) \cdot \bt  + f u \nabla ( f_n u ) \cdot \bt
\eean
Next we want to prove that 
\bean\label{plus}
\lim_{n \ri \infty } \int_{\Omega}  f_n u \nabla ( f u) =
\lim_{n \ri \infty } \int_{\Omega}   f u \nabla ( f_n u )  = \int_{\Omega}   
f^+ u \nabla ( f u).
\eean
The first limit in (\ref{plus}) is clear.
We observe that 
\bea
f_n u \nabla ( f u) - f u \nabla ( f_n u ) = u^2 (f \nabla f_n - f_n \nabla f)
\eea
Since $\Omega$ is two dimensional and $u$ is $H^1_0(\Omega)$, we can assert by the Sobolev embeddings
that $u^2$ is in $L^2(\Omega)$.
$f \nabla f_n - f_n \nabla f$ converges to $f \nabla f^+ - f^+ \nabla f$ in $L^2(\Omega)$, which
is zero, so the second limit in (\ref{plus}) is proved. 
Going back to (\ref{div th}), we have now shown,
$$
\int_{\Omega}  f^+u \nabla ( f u) \cdot \bt =0.
$$
We now multiply equation (\ref{deg 1 eqn}) by $f^+u$, we integrate over $\Omega$, 
and we use that $\alpha$ is non-zero to find that 
$\int_{\Omega}  f^+u^2$ is zero. Similarly, $\int_{\Omega}  f^-u^2$ is zero.
As $f^{-1}(\{ 0 \}) \cap \Omega$ has measure zero, this shows that $u$ is zero. \\
We now consider the case where $\alpha$ is zero. 
After a linear change of variables, we may assume that $\bt$ is the base vector $\bev_1$.
Let $A$  be a  constant such that $\Omega$ is included in the
box $-A \leq x_1, x_2 \leq A$. 
We first note that for any function $g$ in $C^1_c((-A,A) \times (-A,A))$ 
 for any $-A \leq x_1, x_2 \leq A$, 
\bea
|g(x_1, x_2)| \leq (\int_{-A}^A |\p_{x_1} g(x_1,x_2)|^2 dx_1)^{1/2} (2 A)^{1/2},
\eea
thus
\bea
\int_{-A}^A \int_{-A}^A |g(x_1,x_2)|^2 dx_1 dx_2 \leq 4 A^2
\int_{-A}^A \int_{-A}^A |\p_{x_1}  g(x_1,x_2)|^2 dx_1 dx_2,
\eea
and this estimate can be extended to all $g$ in $H^1_0((-A,A) \times (-A,A))$.
It thus follows that $fu$ in zero in $(-A,A) \times (-A,A)$, 
thus $u$ is zero in $\Omega$.

\subsection{Lipschitz stability theorem for a fixed slip}
Let $R$ be a closed rectangle in the plane $x_3=0$.
Let $B$ be a set of triplets $(a, b, d)$ such that the set
\bea 
 \Gamma_{a, b, d} =\{  (x_1, x_2, a x_1 + b x_2 +d): (x_1, x_2) \in R\}
\eea 
is included in the half-space $x_3 <0$.
When appropriate, we will use the short hand notation $m=(a,b,d)$.
We assume that $B$ is a closed and bounded subset of $\RR^3$. It follows that 
that 
\bean \label{pos dis}
\begin{array}{l}
\mbox{\sl the distance between $\Gamma_{m}$ and the plane $x_3=0$ is bounded below}
\\
\mbox{\sl by the same  positive constant for all $ m$ in $B$.}
\end{array}
\eean
We set  $\bn=(-a,-b,1)/\s{1+a^2+b^2}$ to be the normal vector on $\Gamma_m$
and $\sigma = \s{1+a^2+b^2}$ the surface element.
Let $H^1_0(R)^2$ be the space of vector fields $\bg=(g_1,g_2)$  on $R$ with
$H^1_0$ regularity.
Define $\bg_m$ the tangential vector field on $\Gamma_m$, $\bg_m= (g_1,g_2, a g_1 + b g_2)$.
Define the operator
\bean 
A_{m} &:& H^1_0 (R)^2 \ri L^2(V)^3 \no \\
&&\bg \ri \int_R  \bH(\bx, y_1, y_2,  a y_1 + b y_2 +d, \bn)
 \bg_m (y_1, y_2) \sigma d y_1 dy_2 .  \label{Aabd}
\eean
It is clear that $A_{m} $ is linear, continuous, and compact.
Note that due to Theorem \ref{uniq1}, $A_{m}$ is injective. 
In the remainder of this section we fix a non zero $\bh$ in $H^1_0(R)$, and we define
a non-linear function
\bean \label{phi def}
\phi: B \ri L^2(V)^3 \\
\phi(m) = \int_R  \bH(\bx, y_1, y_2,  a y_1 + b y_2 +d, \bn)
 \bh_m (y_1, y_2) \sigma d y_1 dy_2,\label{phi def2}
\eean 
where  $m=(a,b,d)$.
Due to the regularity of the Green's tensor $\bH(\bx, \by, \bn) $, it is clear 
that $\phi$ is real analytic in $m$. 
Now due to Theorem \ref{uniq1}, $\phi$ is injective. 
We now want to prove that the inverse of $\phi$ defined on $\phi(B)$ and valued in $B$
is Lipschitz continuous. This will be achieved by showing that we can apply the inverse function
Theorem.
\begin{thm}\label{stability1}
Fix a non-zero $\bh$ in $H^1_0(R)^2$ and define  the function $\phi$ from $B$ to $L^2(V)^3$
by (\ref{phi def}).
Assume that:\\
(i). $B$ does not contain any triplet in the form $(0,0,d)$, in other words no horizontal 
profiles are allowed for the faults. \\
or\\
(ii). $\bh$ is one- directional.  \\
or\\
(iii). $\bh$ is in $H^2_0(R)^2$.\\
There is a positive constant $C$ such that 
  \bean  
	\label{stabest}
	C |m - m'| \leq   \|  \phi(m) - \phi(m') \|_{L^2(V)} ,
	\eean
	for all $m$ and $m'$ in $B$.
\end{thm}
\textbf{Proof:} Fix $m$ in $B$. Our first task is to evaluate $\nabla \phi (m)$.
We first note that $\bn \sigma$ simplifies to 
 $(-a, -b, 1)$. We recall that $\bH(\bx, \by, \bn)$ is linear in $\bn$.
By the chain rule, for $\by = a y_1 + by_2 +d$,
\bea
\f{\p}{ \p a } \bH(\bx, \by, \bn \sigma) & = &\f{\p y_3}{ \p a } (\p_{y_3}\bH)(\bx, \by, \bn \sigma) -
  \bH(\bx, \by, \bev_1) \\
 & = & y_1 (\p_{y_3}\bH)(\bx, \by, \bn \sigma) -
  \bH(\bx, \by, \bev_1).
\eea
Similarly, 
\bea
 \f{\p}{ \p b } \bH(\bx, \by, \bn) =  y_2 (\p_{y_3}\bH)(\bx, \by, \bn \sigma) -
  \bH(\bx, \by, \bev_2),
\eea and 
\bea
\f{\p}{ \p d } \bH(\bx, \by, \bn) =   (\p_{y_3}\bH)(\bx, \by, \bn \sigma).
\eea
Arguing by contradiction, assume that for some $m$ in $B$, $\nabla \phi (m)$ does not have full rank.
Then there is a non-zero vector $(\gamma_1, \gamma_2, \gamma_3 )$ in 
$\RR^3$ such that 
\bean \label{dependent}
\gamma_1 \f{\p}{ \p a } \phi (m) +  \gamma_2 \f{\p}{ \p b }\phi (m) +  
\gamma_3 \f{\p}{ \p d } \phi (m) =0.
\eean 
Set $f(y_1,y_2) = \gamma_1 y_1  +  \gamma_2 y_2 +  
\gamma_3$.
We note that $ (\gamma_1 \f{\p}{ \p a } +  \gamma_2 \f{\p}{ \p b } +  
\gamma_3 \f{\p}{ \p d }) \bh_m = \nabla f \cdot \bh_m \bev_3$.
Relation (\ref{dependent}) can be expressed as 
\bean \label{onV}
\int_R 
\bH(\bx, y_1, y_2,  a y_1 + b y_2 +d, \bn)
 \nabla f \cdot  \bh_m (y_1, y_2) \bev_3
 \sigma d y_1 dy_2 \no \\
\int_R 
(\p_{y_3}\bH)(\bx, y_1, y_2,  a y_1 + b y_2 +d, \bn)
 \bh_m (y_1, y_2) 
f(y_1,y_2) \sigma d y_1 dy_2 \no \\
- \int_R 
 \bH(\bx, y_1, y_2,  a y_1 + b y_2 +d, \nabla{f})
 \bh_m (y_1, y_2) 
 d y_1 dy_2=0,
\eean
for all $\bx$ in $V$. 
Set $\bw(\bx)$ to be the left hand side of (\ref{onV}) where $\bx$ has been
extended to $\RR^{3-} \setminus \Gamma_{m}$. We now proceed to prove that
$\bw$ is zero in $\RR^{3-} \setminus \Gamma_{m}$.
First, it is clear that $\bw$ satisfies the elasticity equations in 
$\RR^{3-} \setminus \Gamma_{m}$ since the scalar differential operators 
$\p_{y_3}$ and $\p_{x_j}$ commute. Next, due to (\ref{onV}), $\bw$ is zero on $V$.
By construction of Green's tensor $\bH$, for any $\bx$ on the plane $x_3=0$, 
 any $\by$ in $\RR^{3-}$, and any fixed vector $\bp$ in $\RR^3$, 
$$
 T_{\bev_3}(\bx) \bH(\bx, \by, \bp)  =0.
$$
We can thus take a $\p_{y_3}$ derivative and commute the matrix differential operator
$T_{\bev_3}(\bx) $ with the scalar differential operator $\p_{y_3}$ 
to obtain
$$
 T_{\bev_3}(\bx) \p_{y_3} \bH(\bx, \by, \bp)  =0.
$$
It follows that $T_{\bev_3} \bw$ is also zero in $V$ and a Cauchy Kowaleski type argument as in
the proof of Theorem \ref{uniq1}, which was given in
 \cite{volkov2017reconstruction},  shows that $\bw$ must be zero everywhere in 
$\RR^{3-} \setminus \Gamma_{m}$. 
In particular the jump of $\bw$ across $\Gamma_{m}$ must also be zero.
Recall the definition of $\bG$ given by (\ref{Gdef}). 
We note that for any vector $\bv$ in $\RR^3$,
$\bH(\bx,\by,\bv) - \bG(\bx,\by,\bv)$ is smooth for all $\bx$ and
$\by$ in $\RR^{3-}$, see \cite{DV}.
Therefore, the jump across $\Gamma_m$ of
\bean
\int_R 
\bG(\bx, y_1, y_2,  a y_1 + b y_2 +d, \bn)
 \nabla f \cdot  \bh (y_1, y_2) \bev_3
 \sigma d y_1 dy_2 \no \\
\int_R 
(\p_{y_3}\bG)(\bx, y_1, y_2,  a y_1 + b y_2 +d, \bn)
 \bh_m (y_1, y_2) 
f(y_1,y_2) \sigma d y_1 dy_2 \no \\
- \int_R 
 \bG(\bx, y_1, y_2,  a y_1 + b y_2 +d, \nabla{f})
 \bh_m (y_1, y_2) 
 d y_1 dy_2, \label{also zero}
\eean
is also zero. Let us write 
$$ \bev_3 =  \alpha \bn + \bt,$$
where $\bt$ is parallel to $\Gamma_m$.
We now use the fact that the free space Green's function is rotation invariant.
After a change of coordinates by rotation, we can assume that $\Gamma_m$ is horizontal
and $\bt = \beta \bev_1$ (for the sake of lighter notations, the new coordinates will be named in the same way as the old coordinates). 
In the new coordinates we note that $\bh_m \cdot \bev_3 =0$,
and we simply write $\bh$ in place of $\bh_m$.
 The expression
(\ref{also zero}) can be written out as
\bea
\int_{R} \bG (\bx, y_1, y_2, \tilde{d}, \bev_3) (\nabla \tilde{f} \cdot \bh ) (\alpha \bev_3 + \beta \bev_1)
\sigma dy_1 dy_2 \\
+ \alpha  \int_{R} \p_{y_3}\bG (\bx, y_1, y_2, \tilde{d}, \bev_3) \bh  \tilde{f}
\sigma dy_1 dy_2 \\
+ \beta  \int_{R} \p_{y_1}\bG (\bx, y_1, y_2, \tilde{d}, \bev_3) \bh  \tilde{f}
\sigma dy_1 dy_2 \\
-   \int_{R} \bG (\bx, y_1, y_2, \tilde{d}, \nabla \tilde{f}) \bh dy_1 dy_2,
\eea
where $\tilde{f}$ is a non -zero 
affine function.
This must also have a zero jump across $R + \tilde{d}$.
We now proceed to write down the expression for that jump thanks 
Lemma \ref{jump lem2} and formulas (\ref{jump1}-\ref{jump4}) to
find 
\bea
(\nabla \tilde{f} \cdot \bh ) (\alpha \bev_3 + \beta \bev_1) \sigma \\
+ \alpha \sigma \f{\lambda}{ \lambda + 2 \mu } (\di (\tilde{f} \bh)) \bev_3 \\
- \beta \sigma \p_1( \tilde{f} \bh) \\
- \p_3\tilde{f} \bh \\
- \f{\lambda}{ \lambda + 2 \mu } \bh \cdot \nabla \tilde{f}  \bev_3 =0.
\eea
As $\sigma \alpha =1$,  this simplifies  along $\bev_3 $ to 
\bean \label{along e3}
\f{\lambda}{ \lambda + 2 \mu } \tilde{f} \di \bh + \nabla \tilde{f} \cdot \bh =0 .
\eean
 The remaining terms lead to the equation
\bea
\beta \sigma (\nabla \tilde{f} \cdot \bh) \bev_1
- \beta \sigma \p_1 (\tilde{f} \bh) - \p_3 \tilde{f} \bh =0,
\eea
that is, to the system
\bean \label{remaining}
- \beta \sigma \p_1 (\tilde{f}  h_2) -\p_3 \tilde{f} h_2 = 0 \no \\
- \beta \sigma \p_1 (\tilde{f}  h_1) + \beta \sigma  (\p_1 \tilde{f} h_1 +
 \p_2 \tilde{f} h_2) -\p_3 \tilde{f} h_1 = 0
\eean
Assume that condition (i) in the statement of Theorem \ref{stability1}
holds.
Then $\Gamma_m$ is not horizontal, thus
$\beta \neq 0$. Note that 
$\nabla \tilde{f}$ is a constant vector.
Then we can use the first line of equation (\ref{remaining}) in conjunction
to Lemma \ref{deg 1 eqn lem} to find that $h_2 =0$.
Then due to the second line of  (\ref{remaining}) and Lemma \ref{deg 1 eqn lem}, $h_1 =0$.
Thus we showed that $\bh$ is zero
in $H^1_0(R)^2$: contradiction. \\
If condition (ii) in the statement of Theorem \ref{stability1}
holds, 
we set $\bh = u V$, where $u$ is a scalar function and $V$ is a fixed vector
and equation  (\ref{along e3}) simplifies to
$$
\f{\lambda}{\lambda + 2 \mu} V \cdot \nabla (\tilde{f} u) + 
\f{2 \mu}{\lambda + 2 \mu}  (V \cdot \nabla \tilde{f} ) u,
$$
so
equation (\ref{along e3}) in conjunction
to Lemma \ref{deg 1 eqn lem} can be used to show  that $u$
is zero.\\
Now, assume that $\Gamma_m$ is horizontal and that condition (iii) holds. 
In that case $\beta = 0$ and equation 
(\ref{remaining}) is void.
We  also note that here $\alpha = \sigma=1$ and 
 that equation (\ref{along e3}) is still valid.
We use that  the jump of the $\p_{x_3}$ derivative across  $\Gamma_m$ of
\bean
\int_{\Gamma_m} \bG (\bx, y_1, y_2, d, \bev_3) (\nabla f \cdot \bh )  \bev_3 
dy_1 dy_2 \no \\
+  \int_{\Gamma_m} \p_{y_3}\bG (\bx, y_1, y_2, d, \bev_3) \bh  f
 dy_1 dy_2 \no \\
-   \int_{\Gamma_m} \bG (\bx, y_1, y_2, d, \nabla f) \bh dy_1 dy_2
\label{flat}
\eean
is zero.
To complete the proof we apply a change of coordinates by rotation about $\bev_3$ such that
$\nabla f$ becomes parallel to $\bev_1$ in the new coordinates.
By homogeneity, we can then assume that $f(x_1,x_2) = x_1 + \gamma_3$.

We now apply formula (\ref{jump4+1/2}-\ref{jump6}) 
to the (zero) $\p_{x_3}$ jump of
(\ref{flat}) to obtain the following
 equation in the direction of $\bev_1$
\bean
-\p_1(\p_1 f h_1) + \f{3 \lambda + 4 \mu}{ \lambda + 2 \mu} \p_1^2(f h_1) + \p_2^2 (f h_1)
+ 2 \f{\lambda + \mu}{ \lambda + 2 \mu} \p_1 \p_2 (f h_2) \no \\
-\p_1 f ( \f{3 \lambda + 4 \mu}{ \lambda + 2 \mu} \p_1 h_1 + \p_2 h_2) =0. \label{deg2}
\eean
We then eliminate $h_2$ in (\ref{deg2}).
    This is done by using (\ref{along e3}) and observing that as $\p_1 f =1$,
\bea
 \p_1 \p_2 (f h_2) = - \p_1^2 (f h_1) - 2\f{\mu}{\lambda} \p_1 h_1,
\eea
so (\ref{deg2}) reduces to,
as $\ds 1 + 4\f{\mu}{\lambda} \f{\lambda + \mu}{ \lambda + 2 \mu} + 
 \f{3 \lambda + 4 \mu}{ \lambda + 2 \mu}  = 4 + 2 \f{\mu}{\lambda}$,
\bean \label{deg2 becomes}
\p_1^2(f h_1) + \p_2^2 (f h_1) + \p_1 h_1
(- 4 - 2 \f{\mu}{\lambda}) - \p_2 h_2 =0
\eean
We multiply by $f$, use again (\ref{along e3}) and simplify to obtain
 \bean \label{obtained}
f^2 \Delta h_1 + f \p_1 h_1 (-3 - 2 \f{\mu}{\lambda}) + (1 + 2 \f{\mu}{\lambda}) f h_1 =0.
\eean 
Note that this not an elliptic PDE as $f$ may be equal to zero in $\Gamma$.
To show that $h_1$ is zero, fix $\epsilon >0$,  let 
$\Gamma^+ = \{ x \in \Gamma: f(x) > \epsilon \}  $ and  
$\Gamma^- = \{ x \in \Gamma: f(x) < - \epsilon \}  $. 
As $\bh$ is in $H^2_0(\Gamma)$, since $\Gamma^+$ is the intersection 
of $\Gamma$ and a half plane, if it is non empty, the Cauchy Kowaleski Theorem can be applied
to (\ref{obtained}) to claim that $ h_1$ is zero in $\Gamma^+$. 
We carry out the same argument on $\Gamma^-$. Finally we let $\epsilon$ tend to zero: 
this proves that $ h_1$ is zero in $\Gamma$. From there we claim that 
$h_2$ is also zero by recalling (\ref{along e3}) and applying Lemma  \ref{deg 1 eqn lem}.\\

We have thus proved that 
for all $m$ in $B$, $\nabla \phi (m)$  has full rank. 
We now include  the set $B$  in a subset 
$B'$ of $\RR^3$ such that $B'$ is open
and property (\ref{pos dis}) still holds for $B'$.
As for every $m$ in $B'$, 
$\nabla \phi (m)$  has full rank, by the inverse function theorem $\phi$ defines a $C^1$
diffeomorphism from an open neighborhood $U_m$ to its image by $\phi$
on $L^2(V)$.
Thus, there is a positive constant $C_m$ such that for all $m'$ and $m''$ in $U_m$,
\bea 
	C_m |m' - m''| \leq   \|  \phi(m') - \phi(m'') \|_{L^2(V)} .
	\eea
Arguing by contradiction, assume that estimate (\ref{stabest}) fails to be true.
Then there are two sequences $m_n'$ and $m_n''$ in $B$ such that $m_n' \neq m_n''$ and
$\f{\|  \phi(m_n') - \phi(m_n'') \|_{L^2(V)}}{ |m_n'- m_n''|} $ tends to zero.
As $B$ is compact, we may assume after extracting subsequences that 
$m_n'$ converges to $\tilde{m}$ and $m_n''$ converges to $\tilde{\tilde{m}}$.
Since $\phi$ is continuous and injective we must have $\tilde{m} = \tilde{\tilde{m}}$.
But for all $n$ large enough $m_n'$ and $m_n''$ must be in the open 
neighborhood $U_{\tilde{m}}$: contradiction.

\section{Second stability theorem: the case of unknown 
slips} \label{Second stability theorem}

In applications  the slip on $\Gamma$ is unknown, therefore this slip cannot be used
to minimize $  \|  \phi(m) - \phi(m_0) \|_{L^2(V)}$ for $ m$ over $B$ as in (\ref{stabest}) to 
find the geometry $m_0$. Instead,  one has to minimize $\| A_m \bh -  A_{m_0} \bh_0  \|_{L^2(V)}$ 
 over all geometries $m$ and all slips $\bh$. 
The unique minimum is zero and only achieved for $m=m_0$ and $\bh = \bh_0$
according to Theorem \ref{uniq1}.
To obtain Lipschitz stability in $|m-m_0|$ we need to add an additional assumption 
on 
$\bh_0$.
A possible additional assumption is to require that $\bh_0$
 be one directional. Physically, this means that the slip on the fault $\Gamma$ occurs
in only one direction. Interestingly, this condition already appeared in another
theoretical study of destabilization modes of faults, \cite{volkov2010eigenvalue}, as discussed
in section \ref{Fault inverse problem}.\\
Recall the definition (\ref{Aabd}) of operator $A_m$. We will need the following lemma.
\begin{lem}
Let $P_m$ be the orthogonal projection onto $\ov{R(A_m)}$ in $L^2(V)$. 
Fix $m_0$ in $B$. 
Then there is a constant $C$ such that
\bean \label{proj_est}
\| P_m - P_{m_0} \| \leq C |m - m_0|,
\eean
for all $m$ in $B$.
\end{lem}
\textbf{Proof:}
We first note that the closure $\ov{R(A_m)}$ of the range of $A_m$ in $L^2(V)$ is equal to
$\ov{R(A_mA_m^*)}$: this is true because the nullspace $N(A^*_m)$ is equal to $N(A_m A_m^*)$ and we 
can then take the orthogonals of each of this subspace.
If $m$ tends to $m_0$,
it is clear $A_m A^*_m$ is norm convergent to $A_{m_0} A^*_{m_0}$ and that
$\| A_m A^*_m - A_{m_0} A^*_{m_0}\| = O(|m-m_0|)$.
Let ${\cal C}$ be the circle in the complex plane centered at the origin with radius
$\|A_{m_0} A^*_{m_0} \| +1$.
The orthogonal projection on the image of $A_m A^*_m$ can be represented
by the  contour integral as follows, see
 \cite{kato2013perturbation},
\bea
P_m = \f{1}{2 i\pi} \int_{{\cal C}} (   \zeta I - A_m A^*_m)^{-1} d \zeta,
\eea
for all $m$ large enough, and where $I$ is the identity operator in $L^2(V)$.
This leads to (\ref{proj_est}).

\begin{thm}\label{stability2}
Fix a non-zero $\bh_0$ in $H^1_0(R)$ and $m_0$ in $B$.
Assume that $\bh_0$ satisfies one of the two following
additional assumptions:\\
(i). $\bh_0$  is one-directional, that is, $\bh_0$ 
is parallel to a fixed tangential vector.\\
(ii). $\bh_0$  is the gradient of a function $\varphi$ in $H^2(\Gamma)$.\\
Then there exists a positive constant $C$ such that
\bean \label{inf formula}
\inf_{\bh \in H^1_0(R)} \| A_m \bh -  A_{m_0} \bh_0  \|_{L^2(V)}  \geq
 C |m - m_0|,
\eean
for all $m$ in $B$.
\end{thm}
\textbf{Proof:} 
Since $I - P_m$ is an orthogonal projection, 
\bea
\| A_m \bh -  A_{m_0} \bh_0  \|_{L^2(V)}^2 \geq 
\| (I- P_m) ( A_m \bh - A_{m_0} \bh_0)\|_{L^2(V)}^2.
\eea
Since $P_m$ is the orthogonal projection on
$\ov{R(A_m)}$,
  $ P_m A_m \bh =  A_m \bh$, and we obtain
	\bea
\| A_m \bh -  A_{m_0} \bh_0  \|_{L^2(V)}^2 \geq 
\| (I- P_m) A_{m_0} \bh_0)\|_{L^2(V)}^2 = 
\|  P_m  A_{m_0} \bh_0 - A_{m_0} \bh_0\|_{L^2(V)}^2 .
\eea
Arguing by contradiction, assume that there is a sequence $m_n$ in $B$ converging to 
$m_0$ such that 
\bean
\| P_{m_n} A_{m_0} \bh_0 - A_{m_0} \bh_0   \|_{L^2(V)} = o (|m_n - m_0|).
\eean
It clearly follows that
\bea
 \| (I - P_{m_n}) (A_{m_n} - A_{m_0}) \bh_0   \|_{L^2(V)} = o (|m_n - m_0|).
\eea
As
\bea
 \| (P_{m_0} - P_{m_n}) (A_{m_n} - A_{m_0}) \bh_0   \|_{L^2(V)} = o (|m_n - m_0|),
\eea
we may write
\bea
 \| (I - P_{m_0}) (A_{m_n} - A_{m_0}) \bh_0   \|_{L^2(V)} = o (|m_n - m_0|).
\eea
Equivalently,
\bean \label{proj}
 (I - P_{m_0})  \f{(A_{m_n} - A_{m_0})}{|m_n - m_0|} \bh_0  =
 o (1).
\eean
As $\f{m_n - m_0}{|m_n - m_0|}$ is a sequence on the unit sphere of $\RR^3$, after
possibly extracting a subsequence we may assume that it converges to some
$\bq$ with $|\bq|=1$.
Taking the limit as $n \ri \infty$ in (\ref{proj}) we find,
\bea
 (I - P_{m_0}) \p_{\bq} A_{m_0}  \bh_0 = 0,
\eea
thus, there is a $\bg_0$ in $H^1_0(R)$ such that 
\bean \label{with rhs}
\p_{\bq} A_{m_0}   \bh_0 =   A_{m_0}  \bg_0.
\eean
We then set $\bq= (\gamma_1, \gamma_2, \gamma_3)$  as in the proof 
of Theorem \ref{stability1}.
Given the form (\ref{Aabd}) of the operator $A_m$ for $m$ in $B$,
$A_{m_0} \bg_0$ can be extended to a vector field on $\RR^{3-} \setminus \Gamma_{m_0}$
satisfying equations (\ref{uj1}-\ref{uj4}) with $\bg_0$ in place of $\bg$ and 
$\Gamma_{m_0} $ in place of $\Gamma$.
In particular, the normal jump of that extended vector field across $\Gamma_{m_0}$
is zero.
The same argument as in the proof of Theorem \ref{stability1} can then be carried
out to show that $\bh_0$ must satisfy, due to (\ref{with rhs}), a partial differential equation
on $\Gamma_{m_0}$. Due to the $A_{m_0}  \bg_0$ term on the right hand side this equation will be unhelpful 
along any direction which is  tangential to $\Gamma_{m_0}$. However we obtain the same  homogeneous equation
in the normal direction which we write here for $\bh_0$
 \bean \label{along e3 0}
\f{\lambda}{ \lambda + 2 \mu } f \di \bh_0+ \nabla f \cdot \bh_0 =0,
\eean
where this equation was written in a rotated coordinate system such that $\Gamma_{m_0}$ is  parallel to the 
new $x_1, x_2$ plane, $\bh_0$ depends only on the new coordinates $x_1, x_2$
and $f$ is a non-zero affine function whose coefficients depend linearly on 
$\gamma_1, \gamma_2, \gamma_3$. \\
If assumption (i) on $\bh_0$ holds then 
we can apply lemma \ref{deg 1 eqn lem} to claim that $\bh_0 $ is zero: contradiction.\\
If assumption (ii) on $\bh_0$ holds then $\varphi$  satisfies the partial
differential equation
\bean \label{along e3 phi}
\f{\lambda}{ \lambda + 2 \mu } f \Delta  \varphi+ \nabla f\cdot \nabla \varphi =0.
\eean
Let $\mbox{sgn}_0$ be the sign function defined on $\mathbb R$ by: $\mbox{sgn}_0(t)=-1$ 
if $t<0$, $\mbox{sgn}_0(0)=0$ and $\mbox{sgn}_0(t)=1$ if $t>0$.
Multiplying (\ref{along e3 phi}) by $(1+\f{2 \mu}{ \lambda })|f|^{\f{2 \mu}{\lambda}}
\mbox{sgn}_0(f)$, we obtain
\bea
\di \left(|f|^{1+\f{2 \mu}{\lambda}} \nabla \varphi \right)=0.
\eea
As by assumption $\bh_0 = \nabla \varphi$ is in $H^1_0(\Gamma)$,
multiplying  by  $\varphi$ and applying Green's theorem  leads
to
\bean \label{pp}
 \int_\Gamma  |f|^{1+\f{2 \mu}{\lambda}}  |\nabla  \varphi|^2 =0.
\eean
Since $f$ is affine, it vanishes  on a set with low dimensionality. We then deduce from the
 identity \eqref{pp}  that  $\varphi $ is zero: contradiction.

\finproof

\section{Appendix: proof of Lemma \ref{all jumps}}
To show formula (\ref{jump2}), we observe that if $\bx$ is not in $\Gamma$, since
$\bg$ is in  $C^\infty_c(\Gamma)$, integrating by parts we can write
\bea
\int_\Gamma (\p_{y_1}\bG)(\bx, \by, \bn) \bg(\by) \, d y_1 d y_2
= - \int_\Gamma\bG(\bx, \by, \bn)  (\p_{y_1}\bg)(\by) \, d y_1 d y_2,
\eea
and then we can apply formula (\ref{jump1}).
We are not aware of formulas (\ref{jump3}) and (\ref{jump4}) appearing anywhere in the literature, so
we believe that a full proof is called for. 
By a Taylor expansion,
\bean \label{taylor}
\bg(y_1, y_2) = \bg(0,0) + \p_{y_1} \bg (0,0) y_1 + \p_{y_2} \bg (0,0) y_2 + O(\rho^2), 
\eean
where $\rho= \sqrt{ y_1^2 + y_2^2}$. 
Let $R>0$ be small enough so that the circle in the plane $x_3=0$
centered at the origin and with radius $R$ is strictly included in
$\Gamma$.
A long calculation (which we performed thanks to the use of a symbolic
computation software), leads to the following expression 
for  $ \bG(\bx, \by, \bev_1) $ where we 
only indicate twice the odd $x_3$ terms for $x_1=x_2=0$: 
 setting 
$$
A= \f18\,{\frac {\lambda+3\mu}{\pi\,\mu\, \left( \lambda+2\,\mu \right) }}, \quad
B= \f18\,{\frac {\lambda+\mu}{\pi\,\mu\, \left( \lambda+2\,\mu \right) }},
$$
$ \bG(\bx, \by, \bev_1) $ is the product of $(\rho^2 + x_3^2 )^{-5/2}$ and
the matrix whose columns are
\bea
(0, 0, 2\, \left(  \left( {\rho}^{2}+{x_{{3}}}^{2} \right)  \left( A-B
 \right) \lambda-2\,B\mu\, \left( {x_{{3}}}^{2}-2\,{y_{{1}}}^{2}+{y_{{
2}}}^{2} \right)  \right) x_{{3}}
), \\
(0,0,12\,\mu\,By_{{1}}y_{{2}}x_{{3}}),\\
(2\,x_{{3}} \left(  \left( A-B \right) {x_{{3}}}^{2}+{\rho}^{2}A+B
 \left( 5\,{y_{{1}}}^{2}-{y_{{2}}}^{2} \right)  \right) \mu, 12\,\mu\,By_{{1}}y_{{2}}x_{{3}} ,0).
\eea
We note that for $x_3>0$
$$
\int_0^R \f{x_3^3 \rho \, d \rho}{(x_3^2 + \rho^2)^{5/2}} =
\f13\,{\frac { \left( {R}^{2}+{x_{{3}}}^{2} \right) ^{3/2}-{x_{{3}}}^{3
}}{ \left( {R}^{2}+{x_{{3}}}^{2} \right) ^{3/2}}},
$$
 thus
\bean \label{int1}
\lim_{x_3 \ri 0^+} \int_0^{2 \pi}\int_0^R \f{x_3^3 \rho \, d \rho 
\, d \theta 
}{(x_3^2 + \rho^2)^{5/2}} = \f{2 \pi}{3}.
\eean
Similarly
\bean
\lim_{x_3 \ri 0^+} \int_0^{2 \pi}\int_0^R \f{x_3 \rho^3 \cos^2 \theta \, d \rho \, d \theta 
}{(x_3^2 + \rho^2)^{5/2}} = \f{2 \pi}{3},
\eean
and
\bean
\lim_{x_3 \ri 0^+} \int_0^{2 \pi}\int_0^R \f{x_3 \rho^3 \sin^2 \theta   d \rho \, d \theta 
}{(x_3^2 + \rho^2)^{5/2}} = \f{2 \pi}{3},
\eean
while by symmetry
\bean \label{int4}
\lim_{x_3 \ri 0^+} \int_0^{2 \pi}\int_0^R \f{x_3 \rho^3 \sin\theta \cos \theta  d \rho \, d \theta 
}{(x_3^2 + \rho^2)^{5/2}} = 0.
\eean
Thus integrating the matrix $ \bG(\bx, \by, \bev_1) $ times $\bg (0,0)$ over the disk 
in the $x_1$-$x_2$ plane with  radius $R$  centered
at the origin for $x_3>0$ and taking the limit as $x_3$ approaches zero we find, 
\bea 
\left(
\begin{array}{c}
g_3 (0,0) \\
0 \\
 {\frac {\lambda}{\lambda+2\,\mu}} g_1 (0,0) 
\end{array}
\right).
\eea
Given that
\bean \label{int5}
\lim_{x_3 \ri 0^+} \int_0^{2 \pi}\int_0^R \f{x_3^3 \rho^2 \, d \rho 
\, d \theta 
}{(x_3^2 + \rho^2)^{5/2}} = 
\lim_{x_3 \ri 0^+} \int_0^{2 \pi}\int_0^R \f{x_3 \rho^4 \, d \rho 
\, d \theta 
}{(x_3^2 + \rho^2)^{5/2}} = 0,
\eean
using Taylor's expansion (\ref{taylor}), formula (\ref{jump3}) is proved.\\
To prove (\ref{jump4}), 
we perform another calculation aided by the use of symbolic
computation software
to find closed form expressions  for $ \p_{y_3}\bG(\bx, \by, \bev_3) $ where, as previously,
 we 
only indicate twice the odd $x_3$ terms for $x_1=x_2=0$.
It can be written out as  the product of
\bean \label{firstterm}
{\frac {1}{ ( {x_{{3}}}^{2}+\rho^2 ) ^{7/2}}}
\eean
and the three column vectors
\bean \label{twocolj4}
(0, 0, -6\, \left(  \left( A+5\,B \right) {x_{{3}}}^{2}+{\rho}^{2} \left( A-5
\,B \right)  \right) x_{{3}}\mu\,y_{{1}}) \no \\
(0, 0, -6\, \left(  \left( A+5\,B \right) {x_{{3}}}^{2}+{\rho}^{2} \left( A-5
\,B \right)  \right) x_{{3}}y_{{2}}\mu
) \no \\
(-6\,y_{{1}}x_{{3}} \left(  \left( {\rho}^{2}+{x_{{3}}}^{2} \right) 
 \left( A-B \right) \lambda+2\,B \left( -3\,{\rho}^{2}+2\,{x_{{3}}}^{2
} \right) \mu \right)
,  \no \\
-6\,y_{{2}}x_{{3}} \left(  \left( {\rho}^{2}+{x_{{3}}}^{2} \right) 
 \left( A-B \right) \lambda+2\,B \left( -3\,{\rho}^{2}+2\,{x_{{3}}}^{2
} \right) \mu \right)
, 
 0)
\eean
Clearly, by symmetry, 
 \bean \label{int6}
\int_0^{2 \pi}\int_0^R \f{y_j x_3 \rho^3 \, d \rho 
\, d \theta 
}{(x_3^2 + \rho^2)^{7/2}} = 
\int_0^{2 \pi}\int_0^R \f{y_j x_3^3 \rho \, d \rho 
\, d \theta 
}{(x_3^2 + \rho^2)^{7/2}} = 0,
\eean
for $j=1$ or 2. Thus there will be no contribution from $\bg(0,0)$.
Similarly, the cross terms
\bean \label{int7}
\int_0^{2 \pi}\int_0^R \f{y_j y_k x_3 \rho^3 \, d \rho 
\, d \theta 
}{(x_3^2 + \rho^2)^{7/2}} = 
\int_0^{2 \pi}\int_0^R \f{y_j y_k x_3^3 \rho \, d \rho 
\, d \theta 
}{(x_3^2 + \rho^2)^{7/2}} = 0,
\eean
are zero if $j \neq k$.
Now a calculation will show that the following limits hold, 
\bean
\lim_{x_3 \ri 0^+} \int_0^{2 \pi}\int_0^R \f{x_3 \rho^5 \cos^2 \theta \, d \rho \, d \theta 
}{(x_3^2 + \rho^2)^{7/2}} =
\lim_{x_3 \ri 0^+} \int_0^{2 \pi}\int_0^R \f{x_3 \rho^5 \sin^2 \theta \, d \rho \, d \theta 
}{(x_3^2 + \rho^2)^{7/2}} =
 \f{8 \pi }{15}  \no \\
\lim_{x_3 \ri 0^+} \int_0^{2 \pi}\int_0^R \f{x_3^3 \rho^3 \cos^2 \theta \, d \rho \, d \theta 
}{(x_3^2 + \rho^2)^{7/2}} =
\lim_{x_3 \ri 0^+} \int_0^{2 \pi}\int_0^R \f{x_3^3 \rho^3 \sin^2 \theta \, d \rho \, d \theta 
}{(x_3^2 + \rho^2)^{7/2}} =
 \f{2 \pi }{15} \label{4limits}
\eean
We then combine (\ref{firstterm}, \ref{twocolj4}, \ref{4limits}) with Taylor formula (\ref{taylor})
to find a contribution of $\f{\lambda}{\lambda + 2 \mu} 
( \p_{y_1} \bg (0,0)  + \p_{y_2} \bg (0,0) )$ in the direction of $\bev_3$, and
$\nabla_\Gamma g_3 (0,0)$ in the $\bev_1, \bev_2$ plane. 
Higher order terms won't contribute since
\bean
\lim_{x_3 \ri 0^+} \int_0^R \f{x_3 \rho^6 \, d \rho \,
}{(x_3^2 + \rho^2)^{7/2}} \, =
 \, \lim_{x_3 \ri 0^+} \int_0^R \f{x_3^3 \rho^4  \, d \rho \,
}{(x_3^2 + \rho^2)^{7/2}} =0,
\eean
and formula (\ref{jump4}) is proved. \\
Formula (\ref{jump4+1/2}) is derived likewise. \\
To prove formula (\ref{jump5}) 
we need  a higher order Taylor formula.  We write, 
\bean \label{taylor2}
\bg(y_1, y_2) = \bg(0,0) + \p_{y_1} \bg (0,0) y_1 + \p_{y_2} \bg (0,0) y_2 \no
\\ + 
\f12 \p_{y_1}^2 \bg (0,0) y_1^2 + \f12 \p_{y_2}^2 \bg (0,0) y_2^2 + 
 \p_{y_1} \p_{y_2} \bg (0,0) y_1 y_2 + 
O(\rho^3).
\eean
We perform another calculation aided by the use of symbolic
computation software
to find expressions for 
$ \p_{x_3} \p_{y_3} \bG(\bx, \by, \bev_3) $.
For the sake of brevity we only give a proof in the case where $g_3 =0 $, 
so we only need  the first two columns.
 We 
only indicate twice the odd $x_3$ terms for $x_1=x_2=0$.
They are the  the product of
$ \ds
\f{1}{(\rho^2 + x_3^2)^{\f92}}
$
and a matrix whose
first  column is
\bea
6\,\mu \left( 2\,{x_{{3}}}^{4}- \left( 41\,{y_{{1}}}^{2}+{y_{{2}}}^{2}
 \right) {x_{{3}}}^{2}+3\,{\rho}^{2} \left( 9\,{y_{{1}}}^{2}-{y_{{2}}}
^{2} \right)  \right) x_{{3}}B+6\, \left( 3\,{\rho}^{4}+{\rho}^{2}{x_{
{3}}}^{2}-2\,{x_{{3}}}^{4} \right) \mu\,x_{{3}}A
\\
-60\,y_{{1}}y_{{2}} \left( -3\,{\rho}^{2}+4\,{x_{{3}}}^{2} \right) Bx_
{{3}}\mu
\\
0,
\eea
and whose second column is
\bea
-60\,y_{{1}}y_{{2}} \left( -3\,{\rho}^{2}+4\,{x_{{3}}}^{2} \right) Bx_
{{3}}\mu
 \\
6\,\mu \left( 2\,{x_{{3}}}^{4}- \left( 41\,{y_{{1}}}^{2}+{y_{{2}}}^{2}
 \right) {x_{{3}}}^{2}+3\,{\rho}^{2} \left( 9\,{y_{{1}}}^{2}-{y_{{2}}}
^{2} \right)  \right) x_{{3}}B+6\, \left( 3\,{\rho}^{4}+{\rho}^{2}{x_{
{3}}}^{2}-2\,{x_{{3}}}^{4} \right) \mu\,x_{{3}}A
\\
0.
\eea
Further calculations show that if we multiply these two columns
by $\ds\f{ \rho}{(\rho^2 + x_3^2)^{\f92}} \ds $ integrate in $\rho$ from 0 to $R$ and 
and $\theta$ from 0 to $2 \pi$
and then take the limit as $x_3$ tends to zero, we find zero.
We also find zero as we multiply by $\ds \f{ \rho^2 \cos \theta}{(\rho^2 + x_3^2)^{\f92}}  $,
$\ds \f{ \rho^2 \sin \theta}{(\rho^2 + x_3^2)^{\f92}}  $, or by 
$\ds \f{ \rho^p}{(\rho^2 + x_3^2)^{\f92}}  $, if $p \geq 4$.
 We only find non-zero contributions (which are independent of $R>0$) for the second derivatives of $\bg$: they are found by
multiplying by $\ds \f{ \rho^3 \sin^p \theta \cos^q \theta}{(\rho^2 + x_3^2)^{\f92}}  $,
where $p,q$ are in $\{0, 1, 2\}$ with $p+q=2$. More precisely, 
for the $\p_1^2$ derivative, we find, for the first column
\bea
({\frac {6\,\lambda+8\,\mu}{\lambda+2\,\mu}}, 0,0),
\eea
for the second column
\bea
(0, 2,0).
\eea
For the $\p_2^2$  derivative, we find, for the first column
\bea
(2, 0,0),
\eea
for the second column
\bea
(0, {\frac {6\,\lambda+8\,\mu}{\lambda+2\,\mu}},0).
\eea
For the  $\p_1 \p_2$ derivative, we find, for the first column
\bea
(0, {\frac {2\,\lambda+2\,\mu}{\lambda+2\,\mu}},0),
\eea
for the second column
\bea
({\frac {2\,\lambda+2\,\mu}{\lambda+2\,\mu}}, 0,0).
\eea
Finally, we arrive at jump formula (\ref{jump5}) thanks to a linear combination. \\
The same proof technique is used for deriving formula (\ref{jump6}). 
For the sake of brevity,  here too we assume $g_3 =0$.
 This time a Taylor expansion of order 1 is sufficient. The first two columns of $ \p_{x_3}\bG(\bx, \by, \bev_1) $ where we 
only indicate twice the odd $x_3$ terms for $x_1=x_2=0$ are
is the product of
$ \ds
\frac {1}{ ({x_{{3}}}^{2}+\rho^2 ) ^{7/2}}
$
and
\bean \label{col1}
6\,x_{{3}}y_{{1}} \left(  \left( {\rho}^{2}+{x_{{3}}}^{2} \right) 
 \left( A-B \right) \lambda+\mu\, \left( 2\,A \left( {\rho}^{2}+{x_{{3
}}}^{2} \right) -2\,B \left( 2\,{x_{{3}}}^{2}-3\,{y_{{1}}}^{2}+2\,{y_{
{2}}}^{2} \right)  \right)  \right) \\
6\,x_{{3}}y_{{2}} \left(  \left( {\rho}^{2}+{x_{{3}}}^{2} \right) 
 \left( A-B \right) \lambda-2\,B\mu\, \left( {x_{{3}}}^{2}-4\,{y_{{1}}
}^{2}+{y_{{2}}}^{2} \right)  \right)\\
0,
\eean
and
\bean \label{col2}
6\,\mu\,x_{{3}}  y_{{2}} \left(  \left( A-B \right) {y_{{2}}}^{2}+ \left( A+9\,
B \right) {y_{{1}}}^{2}+{x_{{3}}}^{2} \left( A-B \right)  \right) \\
6 \mu x_{{3}} \,y_{{1}} \left( {x_{{3}}}^{2} \left( A-B \right) + \left( A-B \right) {y
_{{1}}}^{2}+{y_{{2}}}^{2} \left( A+9\,B \right)  \right) \\
0.
\eean
We then  multiply (\ref{col1}) and (\ref{col2}) by $\ds\f{ \rho y_i}{(\rho^2 + x_3^2)^{\f72}}$, $i=1,2$ and integrate in $\rho$ from 0 to $R$, 
and $\theta$ from 0 to $2 \pi$ to only find non-zero contributions (which are independent of $R>0$) for the first  derivatives of $\bg$.
After simplification, we obtain formula (\ref{jump6}).

\section{Conclusion}
In this paper we studied the well-posedness of the  fault inverse problem.  We derived 
 stability  estimates  for determining the plane containing the fault. 
We proved that 
 if the  slip field is  known,  this determination is Lipschitz stable.  
 In the more realistic  case where the slip field is unknown,  we  showed  another
Lipschitz stability result  under the additional assumption, which seems physically relevant,
that the slip field is one directional.  The proofs of the  results presented in this paper
 are non-constructive 
and thus provide no insight on how
 the  stability constants depend on  the physics and on the geometry
 of the problem. This will be the subject of forthcoming work.

\section*{Acknowledgements}
D. Volkov is supported  by
a Simons Foundation Collaboration Grant. The research of FT was supported in part by the
 grant ANR-17-CE40-0029 of the French National Research Agency ANR (project MultiOnde),
 and the LabEx PERSYVAL-Lab (ANR-11-LABX- 0025-01).


\end{document}